\documentclass[11pt]{amsart}

\def\grade{\operatorname{grade}}
\def\depth{\operatorname{depth}}
\def\dim{\operatorname{dim}}
\def\ZZ{{\mathbb Z}}

\def\mm{{\mathfrak m}}
\def\lrar{{\longrightarrow}}

\def\G{{\mathcal G}}
\def\E{{\mathcal E}}
\def\F{{\mathcal F}}
\def\lm{{\lambda}}
\def\bl{\begin{lemma}}
\def\el{\end{lemma}}
\def\bt{\begin{theorem}}
\def\et{\end{theorem}}
\def\ben{\begin{enumerate}}
\def\een{\end{enumerate}}
\def\bpf{\begin{proof}}
\def\epf{\end{proof}}
\def\beqn{\begin{eqnarray*}}
\def\eeqn{\end{eqnarray*}}
\def\bd{\begin{definition}}
\def\ed{\end{definition}}
\def\bp{\begin{proposition}}
\def\ep{\end{proposition}}
\def\bc{\begin{corollary}}
\def\ec{\end{corollary}}

\newtheorem{lemma}{Lemma}
\newtheorem{corollary}[lemma]{Corollary}
\newtheorem{theorem}[lemma]{Theorem}
\newtheorem{proposition}[lemma]{Proposition}

\newtheorem{definition}[lemma]{Definition}

\advance\headheight1.15pt

\begin{document}

\title {Hilbert coefficients and depths of form rings}

\author[A. V. Jayanthan]{A. V. Jayanthan$^*$}
\author{Balwant Singh \and J. K. Verma}
\thanks{$^*$Supported by the National Board for Higher
Mathematics, India.}
\address{Department of Mathematics, Indian Institute of Technology 
Bombay, Powai, Mumbai, India - 400076}
\email{jayan@math.iitb.ac.in}
\email{balwant@math.iitb.ac.in}
\email{jkv@math.iitb.ac.in}

\maketitle

\begin{abstract}
We present short and elementary proofs of two theorems of Huckaba and
Marley, while generalizing them at the same time to the case of a
module. The theorems concern a characterization of the depth of the
associated graded ring of a Cohen-Macaulay module, with respect to a
Hilbert filtration, in terms of the Hilbert coefficient $e_1$. As an
application, we derive bounds on the higher Hilbert coefficient $e_i$
in terms of $e_0.$
\end{abstract}

\maketitle

\section{Introduction}
Let $(R, \mm)$ be a Noetherian local ring with infinite residue field
and let $M$ be a finitely generated $R$-module.  Let $\F = \{I_n\}_{n 
\geq 0}$ be
a filtration of ideals of $R$. The filtration $\F = \{I_n\}_{n \geq 0}$
is called a {\it
Hilbert filtration with respect to $M$} if it is $I_1$-good with
respect to $M$ (i.e. $I_1I_nM = I_{n+1}M$ for large $n$) and $I_1$ is
an ideal of definition for $M$. Let $G(\F,M) =
\bigoplus_{n \geq 0}I_nM/I_{n+1}M$ be the associated graded module of
$M$ with respect to $\F.$ In \cite{h}, Huckaba proved a
characterization for $G(\F,M)$ to have depth at least $d-1,$ for the
$I$-adic filtration and $M = R$. The characterization is given in
terms of $e_1(I),$ the first Hilbert coefficient of $R$ with respect
to $I.$ (See Section 1 for the definition of $e_1$ and other Hilbert
coefficients.) Huckaba's proof is based on a generalized version of
the ``Fundamental Lemma" of Huneke.  In \cite{hm}, Huckaba and Marley
used a modification of the Koszul complex to prove the following
general version of Huckaba's result and also a characterization for 
$G(\F) = G(\F, R)$ to be Cohen-Macaulay. 

\noindent
{\bf Theorem A. (Huckaba-Marley)} Let $(R,\mm)$ be a  
Cohen-Macaulay local ring of dimension $d$ with infinite 
residue field and $d\geq 1.$ 
Let $\F = \{I_n\}_{n \geq 0}$ be a Hilbert filtration in $R$, and let
$J$ a minimal reduction of $\F$.  Then 
\begin{enumerate}
\item[(i)] $e_1(\F) \geq \sum_{i \geq 1}\lm(I_i + J/J)$ and the 
equality holds if and only if $G(\F)$ is Cohen-Macaulay.
\item[(ii)] $e_1(\F) \leq \sum_{i \geq 1}\lm(I_i/JI_{i-1})$ and the
equality holds if and only if $\depth G(\F) \geq d-1$.
\end{enumerate}
In this article we present, among other things, short and elementary 
proofs of both the results in a slightly more general setting, namely 
for a Cohen-Macaulay module  over a Noetherian local ring. Our key 
observation is that the Hilbert 
coefficient $e_1$ of $M$ with respect to $\F$ can be expressed as the 
Hilbert coefficient $e_0$ of a suitable module over the Rees ring of 
$R$ with respect to $\F.$ This leads to a simpler proof of the main 
results, as the coefficient $e_0$ is much better understood than the 
higher coefficients.

\section{The module $\E_J(\F,M)$}

Let $(R, \mm)$ be a Noetherian local ring with infinite residue field
and let $M$ be a finitely generated $R$-module of dimension $d\geq 1$.
Let $\F = \{I_n\}_{n \geq 0}$ be a Hilbert filtration with respect to
$M$. Let $J \subseteq I_1$ be an ideal of
$R$. We say that $J$ is a reduction of $\F$ with respect to $M$ if
there exists an integer $r$ such that $JI_nM = I_{n+1}M$ for all $n 
\geq r$. A reduction is called a minimal reduction if it is minimal
with respect to inclusion.
The Hilbert coefficients $e_i(\F,M)$ of $M$ with respect to the 
Hilbert filtration $\F$ are defined by writing the 
Hilbert-Samuel polynomial $P_{\F}(M,n)$ corresponding to the Hilbert function 
$n\mapsto H_{\F}(M,n) = \lambda(M/I_nM)$, where $\lambda$ denote length 
as $R$-module in the following form : 
$$P_\F(M,n) = \sum_{i=0}^{d}(-1)^ie_i(\F,M){n+d-1-i \choose d-i}. $$ 

We write $e_i(\F)$ for  $e_i(\F,R).$ 
In the following lemma we prove the existence of a minimal
reduction for $\F$ with respect to $M$. For lack of a suitable
reference in the literature, we include a proof.

\begin{lemma}
Let $(R,\mm)$ be a Noetherian local ring with infinite residue field
and let $M$ be a finite $R$-module of dimension $d$. Let $\F$ be a
Hilbert filtration with respect to $M$. Then there exist $x_1, \ldots,
x_d \in I_1$ such that $({\bf \underline{x}}) = (x_1, \ldots, x_d)$ is
a minimal reduction of $\F$ with respect to $M$ and $e_0(({\bf
\underline{x}}), M) = e_0(\F, M)$.
\end{lemma}
\begin{proof}
By Corollary 4.6.10 of \cite{bh}, there exist $x_1, \ldots, x_d \in
I_1$ such that $({\bf\underline{x}})=(x_1, \ldots, x_d) $ is a minimal reduction of
$I_1$ with respect to $M$ and $e_0(({\bf\underline{x}}), M) =
e_0(I_1, M)$. Since $\F$ is $I_1$-good, there exists an integer $g$
such that $I_1^k I_n = I_{n+k}$ for all $k \geq 1$ and $n \geq g$.
Choose $n \gg 0$ and $k \geq r$, where $r$ is the minimal integer $n$
which satisfies the equation $({\bf\underline{x}})I_1^nM =
I_1^{n+1}M$. Then $I_{n+k+1}M = I_1^{k+1}I_nM
=({\bf\underline{x}})I_1^kI_nM = ({\bf\underline{x}})I_{n+k}M$. Thus
$({\bf\underline{x}})$ forms a minimal reduction of $\F$ with respect
to $M$. Since $\F$ is an $I_1$-good filtration, the second statement
also follows. 
\end{proof}

\vskip 3mm
Let $J$ be a $d$-generated 
minimal reduction of $\F$ with respect to $M$ and let $r$ denote the
reduction number of $\F$ with respect to $J$ and $M$, i.e. $r$ is the
smallest integer $n$ such that $JI_nM=I_{n+1}M.$ 

Recall that the extended Rees algebra of $\F$ is the graded algebra 
$\E(\F)=\bigoplus _{i\in \ZZ} I_it^i,$ where $t$ is an indeterminate 
and $I_i=R$ for $i\leq 0.$ The extended Rees algebra  $\E(J)$ of $J$ 
is a graded subalgebra of $\E(\F).$ For a finitely generated $R$-module
$M$, let $\E(\F, M)$ denote the module $\bigoplus_{i \in \ZZ}I_iMt^i$.
Then $\E(J,M)$ (:= $\E(\G, M)$, where $\G$ is the $J$-adic filtration)
is a graded submodule of $\E(\F, M)$. Also,
$\E(\F,M)$ is a finitely generated graded module over $\E(J)$. Put
$\E_J(\F,M) = \E(\F,M)/\E(J,M)$. Write 
$$\E_J(\F,M) = \bigoplus_{i \in \ZZ} 
\E_J(\F,M)_it^i \;\;\mbox{with}\;\; \E_J(\F,M)_i = I_iM/J^iM.$$ 
Note that $\E_J(\F,M)$ is a finitely generated graded
$\E(J)$-module and that $\E_J(\F,M)_i
= 0$ for $i \leq 0.$ 
Let $\gamma(\F,M)=\grade (G(\F)_+, G(\F,M)).$ It is well known that 
$\gamma(\F,M)=\depth _{G(\F)} (G(\F,M)).$ We write 
$\gamma(\F)$ for $\gamma(\F,R).$ 
For the $\E(J)$-module $\E_J(\F,M),$ let 
$$e_i(\E_J(\F,M)) \;\mbox{ (respectively}\;\; e_i(Jt,\E_J(\F,M)))$$ 
denote the Hilbert coefficients of the Hilbert polynomial associated 
to the Hilbert function 
$$n\mapsto \lm(\E_J(\F,M)_n)\;\; \mbox{(respectively}\; n\mapsto 
\lm(\E_J(\F,M)/(Jt)^n\E_J(\F,M))).$$

\bp \label{p1}Let $M$ be Cohen-Macaulay and  
$\E_J(\F,M) \neq 0$. Then 
\begin{enumerate}
\item [(1)] $\lm(\E_J(\F,M)/Jt\;\E_J(\F,M)) = \sum_{n=1}^{\infty}
\lm(I_nM/JI_{n-1}M)
< \infty.$

\item [(2)] $\dim \E_J(\F,M) = d$ and $e_i(\F,M)=e_{i-1}(\E_J(\F,M))$ 
for every 
$i,\; 1\leq i\leq d.$ In particular, $e_1(\F, M)> 0.$  

\item [(3)] $e_1(\F,M)=e_0(\E_J(\F,M))=e_0(Jt,\E_J(\F,M)).$ 

\item [(4)] $\gamma(\F,M)\geq \depth \E_J(\F,M) -1.$
\end{enumerate}
\ep

\begin{proof}(1) The equality is clear from the definition of 
$\E_J(\F,M),$ and the sum is finite because $I_nM=JI_{n-1}M$ for $n \gg 0.$  

(2) Choose least $r$ such that $I_rM \neq J^rM.$ Let 
$c \in I_rM \setminus J^rM.$ Then by minimality of $r,$ we have 
$c \in J^{r-1}M.$ Let $ \bar c$ denote the image of $c$ in 
$J^{r-1}M/J^rM.$ Since $J$ is generated by 
an $M$-regular sequence of length $d,\; G(J,M) \cong M/JM[X_1, \ldots, 
X_d].$ Therefore for $n \geq r,$ the $R/J$-submodule 
$(cJ^{n-r}M+J^nM)/J^nM$ of 
$J^{n-1}M/J^nM$ is generated minimally by the set $\{\bar 
cX^{\alpha}\mid \;\mid\alpha\mid=n-r \},$ whose cardinality is 
${n+d-r-1\choose d-1}.$ Therefore 
$$\lm ((cJ^{n-r}M+J^nM)/J^nM)\geq {n+d-r-1\choose d-1}.$$ Since 
$ cJ^{n-r}M+J^nM \subseteq I_rJ^{n-r}M\subseteq I_nM,$ we get 
$$\lm(\E_J(\F,M)_n)=\lm (I_nM/J^nM)\geq {n+d-r-1\choose d-1}.$$ Hence
$\dim \E_J(\F,M) \geq d.$ 
Since $G(J, M) \cong M/JM[X_1, \ldots, X_d]$,  
$$\lm(M/J^nM )= \lm(M/JM) {n+d-1 \choose d} =  e_0(\F,M) {n+d-1 
\choose d}$$ for every $n.$ 
On the other hand, for large $n$ we have
$$\lm(M/I_nM )= \sum_{k=0}^d(-1)^k e_k(\F,M) {n+d-1-k \choose d-k},\; \mbox{and} $$
$$\lm(\E_J(\F,M)_n)= \sum_{k=0}^{d-1}(-1)^k e_{k+1}(\F,M) {n+d-2-k \choose 
d-1-k}.$$ 
This shows that $\dim (\E_J(\F,M))\leq d.$ Combining this with 
the inequality $\dim (\E_J(\F,M))\geq d,$ all the assertions of (2) follow.

(3) Let $\Gamma_n=\bigoplus _{i\geq n}\E_J(\F,M)_it^i.$ Then 
$(Jt)^n\E_J(\F,M) \subseteq \Gamma_{n+1}.$ On the other hand, if 
$\E_J(\F,M)$ is generated as an $\E(J)$-module by 
$$\E_J(\F,M)_1t\bigoplus \cdots \bigoplus \E_J(\F,M)_st^s$$ then 
$\Gamma_n\subseteq 
(Jt)^{n-s}\E_J(\F,M)$ for $n\geq s.$ Thus 
for all large $n$ we have $\Gamma_n\subseteq (Jt)^{n-s}\E_J(\F,M)\subseteq
\Gamma_{n-s+1},$ whence 
\begin{eqnarray*}
\lm(\E_J(\F,M)/\Gamma_n)
&\geq &\lm(\E_J(\F,M)/(Jt)^{n-s}\E_J(\F,M))\\
&\geq& \lm(\E_J(\F,M)/\Gamma_{n-s+1}).
\end{eqnarray*}
 It 
follows that $e_0(\E_J(\F,M))=e_0(Jt,\E_J(\F,M)).$ We have already proved 
the equality 
$e_1(\F,M)=e_0(\E_J(\F,M))$  in (2).

(4) This follows from Proposition (1.2.9) of \cite{bh} in view of the 
exact sequence 
$$0 \lrar \E(J,M) \lrar \E(\F,M) \lrar  \E_J(\F,M) \lrar 0$$ 
of $\E(J)$-modules. 
\end{proof}

\section{Main Theorems}

\vspace*{1mm} 

We use the properties of $\E_J(\F,M)$ proved in Section 1 to prove 
the following generalized version of part (ii) of Theorem A.

\vspace{3mm} 

\begin{theorem}\label{d-1}
Suppose $M$ is a finitely generated Cohen-Macaulay $R$-module of dimension
$d\geq 1$ and let $\F$ a Hilbert filtration with respect to $M$.
Then $e_1(\F,M) \leq \sum_{n=1}^{\infty} \lambda(I_nM/JI_{n-1}M)$ for
every minimal reduction $J$ of $\F.$ Moreover, the following three
conditions are equivalent:
\begin{enumerate}
\item [(1)] $e_1(\F,M) = \sum_{n=1}^{\infty}\lambda(I_nM/JI_{n-1}M)$ for
some 
$($resp. every$)$ minimal 
reduction $J$ of $\F.$
\item [(2)] $\E_J(\F,M)$ is Cohen-Macaulay for some $($resp. every$)$ minimal 
reduction $J$ of $\F.$
 
\item [(3)] $\gamma(\F,M) \geq d-1.$
\end{enumerate}
\end{theorem}

\begin{proof}
If $\E_J(\F,M) = 0$ then the inequality and the three conditions hold. Assume
therefore that $\E_J(\F,M) \neq 0.$  
Then by Proposition \ref{p1} we have $e_1(\F,M) = e_0(Jt, \E_J(\F,M)).$ 
Since $Jt$ is clearly generated by a system of parameters for 
the $\E(J)$-module $\E_J(\F,M),$ we have 
$e_0(Jt, \E_J(\F,M))\leq \lambda(\E_J(\F,M)/Jt\E_J(\F,M))$ 
and the last quantity equals 
$\sum_{n=1}^{\infty}\lambda(I_nM/JI_{n-1}M)$ by Proposition  \ref{p1}.
This proves the inequality of the theorem.

Now, $\E_J(\F,M)$ is Cohen-Macaulay if and only if $Jt$ is 
generated by a
regular sequence on $\E_J(\F,M)$ if and only if $e_0(Jt, \E_J(\F,M)) =
\lambda(\E_J(\F,M)/Jt\E_J(\F,M)).$ By Proposition \ref{p1} The last 
condition is equivalent to 
$e_1(\F,M) = \sum_{n=1}^{\infty}\lambda(I_nM/JI_{n-1}M).$ This proves the 
equivalence of (1) and (2). 
The equivalence of (2) and (3) follows from Proposition 1.2.9 of
\cite{bh} in view of the exact sequence 
$$
0 \lrar \E(J,M) \lrar \E(\F,M) \lrar \E_J(\F,M) \lrar 0 
$$
of $\E(J)$-modules.
\end{proof}

\vspace{3mm} 

Now, in order to prove part (ii) of the results of Huckaba-Marley in 
the case of a module we need the  following lemma for the induction 
procedure.
This lemma is a generalization for modules of a  special case of 
Lemma 2.2 of \cite{hm}. For a nonzero element $x$ of $R$ let $x^*$ 
denote its initial form in $G(\F)$.

\vspace{3mm} 

\begin{lemma}\label{nzd}
Let $M$ be a finitely generated  
$R$-module, and let $\F$ be a Hilbert filtration with respect to $M$.
Let $x \in I_1 \backslash I_2$ be superficial for $\F$ and $M$. If
$\gamma(\F/xR, M/xM) > 0$, then $x^* \in G(\F)$ is regular 
on $G(\F,M).$
\end{lemma} 

\begin{proof}
Let $y\in I_t\backslash I_{t+1}$ be such that $(\bar{y})^* \in G(\F/xR)$ 
is regular on $G(\F/xR, M/xM).$ Then $(I_{n+tj}M :_M y^j) 
\subset (I_n, x)M$ for
all $n, j$. Since $x$ is superficial for $\F$ and $M$, there exists an 
integer
$c$ such that $(I_{n+j}M :_M x^j) \cap I_cM = I_nM$ for all $n > c, \, 
j \geq
1$. Let $n$ and $j$ be arbitrary and let $p > c/t$. Then 
$$
y^p(I_{n+j}M :_M
x^j) \subseteq (I_{n+pt+j}M :_M x^j) \cap I_cM \subseteq I_{n+pt}M.
$$ 
Therefore 
$(I_{n+j}M :_M x^j) \subseteq I_{n+pt}M :_M y^p \subseteq (I_n, x)M$. Thus
$$I_{n+j}M :_M x^j = I_nM + x(I_{n+j}M :_M x^{j+1}).$$ Iterating this formula
$n$-times, we get 
$$
I_{n+j}M :_M x^j = I_nM + xI_{n-1}M + \cdots + x^n(I_{n+j}M :_M x^{j+n})
= I_nM .$$
Therefore $x^*$ is regular on $G(\F,M)$.
\end{proof} 

\vspace*{3mm}

\begin{theorem}\label{d}
Let $M$ be a finitely generated Cohen-Macaulay $R$-module of dimension
$d\geq 1$. Let $\F$ be a Hilbert filtration with respect to $M$
and let $J \subseteq I_1$ be a minimal reduction of $\F$
with respect to $M$ with reduction number $r$.  Then 
\noindent $e_1(\F,M) \geq \sum_{i = 1}^r\lm(I_iM + JM/JM),$ and 
the equality holds if
and only if $G(\F,M)$ is Cohen-Macaulay as a $G(\F)$-module.
\end{theorem}

\begin{proof} 
We use induction on $d.$ Let $d=1$ and $J=aR.$
We set 
$$\alpha_i=\lambda(I_iM/J^iM),\; 
\beta_i =\lambda((I_iM\cap JM)/JI_{i-1}M),$$ 
and $\gamma_i=\lambda(I_iM/(I_iM\cap JM)).$
Note that $\beta_i=\gamma_i=0$ 
for $i>r.$ We claim that 
$\alpha_n= \sum_{i=1}^n (\beta_i+\gamma_i)$  for all $n\geq 1.$ This 
holds trivially for 
$n=1.$ We have $\beta_n+\gamma_n= \lambda(I_nM/JI_{n-1}M).$ Therefore 
the exact sequence 
$$
0\rightarrow I_{n-1}M/J^{n-1}M \stackrel{a}{\rightarrow} 
I_nM/J^nM\rightarrow I_nM/JI_{n-1}M \rightarrow 0 
$$ 
 gives 
$\alpha_n= \beta_n+\gamma_n + \lambda(I_{n-1}M/J^{n-1}M),$ proving the 
claim by induction on $n.$ 
Since $\dim (\E_J(\F,M))=1$,  we have
$$e_1(\F,M)=e_0(\E_J(\F,M))=\lambda(\E_J(\F,M)_n)$$ for $n \gg 0$. Hence
for large $n$,
$$e_1(\F,M)=\lambda(I_nM/J^nM)=
 \sum_{i=1}^n (\beta_i+\gamma_i) =
 \sum_{i=1}^r (\beta_i+\gamma_i).$$
 Thus $e_1(\F,M)\geq \sum_{i=1}^r 
\gamma_i$ with equality if and only if each $\beta_i=0$ if 
and only if $a^{*}\in I_1/I_2$ is $G(\F,M)$-regular. This 
proves the assertion for $d=1.$ 

Assume now that $d > 1$ and that $G(\F,M)$ is Cohen-Macaulay.
Let $J = (a_1, \ldots, a_d)$ be a reduction of $\F$ with respect to M
such that $a_1^*, \ldots, a_d^* \in G(\F)$ is a regular sequence on
$G(\F, M)$. Put  $\bar{\F}=\F/a_1R $ and $\bar{M} = 
 M/a_1M.$ 
  Then $G(\F, M)/a_1^*G(\F, M)  \cong G(\bar{\F}, \bar{M}).$ 

Hence $G(\bar{\F}, \bar{M})$ is Cohen-Macaulay. By induction
\begin{eqnarray*}
e_1(\bar{\F}, \bar{M}) & = & \sum_{i \geq 1} \lm(\bar{I}_i\bar{M} + 
\bar{J}\bar{M}/\bar{J}\bar{M}) \\
&  = & \sum_{i \geq 1}\lm(I_iM+JM+a_1M/JM +a_1M)\\ 
& = & \sum_{i \geq 1} \lm(I_iM + JM/JM).
\end{eqnarray*}
Since $a_1^*$ is regular on $G(\F,M)$, we have $e_1(\F,M) =
e_1(\bar{\F}, \bar{M}),$ and we get the desired equality.  

Conversely, let $d > 1$ and let $e_1(\F,M)$ have the given form.
Choose $a_1,\ldots, a_d \in I_1 $ such that
$a_1$ is superficial for $\F$ and $M$ and $J = (a_1, \ldots, a_d)$. Then

\begin{eqnarray*}
e_1(\F,M) &=& e_1(\bar{\F}, \bar{M})\\
          & \geq& \sum_{n \geq 1}
\lm\left(\frac{\bar{I}_n\bar{M}+ \bar{J}\bar{M}}{\bar{J}\bar{M}}\right)
 =   \sum_{n \geq 1} \lm(I_nM+JM/JM). 
\end{eqnarray*}
Now, suppose that the equality holds. Then

\begin{eqnarray*}
e_1(\bar{\F}, \bar{M})& =& e_1(\F, M)\\
& =& \sum_{n \geq 1} \lm(I_nM+JM/JM) 
=  \sum_{n \geq 1}\lm\left(\frac{\bar{I}_n\bar{M}+ 
\bar{J}\bar{M}}{\bar{J}\bar{M}}\right).
\end{eqnarray*}

By induction $G(\bar{\F}, \bar{M})$ is Cohen-Macaulay. Thus
$\gamma(\bar{\F}, \bar{M})> 0$. Therefore by Lemma 
\ref{nzd}, $a_1^* \in G(\F)$ is $G(\F, M)$-regular  and hence 
$G(\F, M)$ is Cohen-Macaulay.
\end{proof} 

\noindent
{\bf Remark :} Huckaba and Marley have derived several consequences of
Theorem A in their paper \cite{hm}. These results can be formulated
for modules in an obvious way.

\section{An Application}

\vspace{1mm} 

\begin{theorem}
Let $(R,\mm)$ be a $d$-dimensional Cohen-Macaulay local ring. Let $I$ be
an $\mm$-primary ideal of $R$ such that $\gamma(I)\geq d-1$. Let
$e_0, e_1 \ldots, e_d$ be the Hilbert coefficients of $R$ with 
respect to $I$. Then 
$$
e_i \leq {e_0 \choose 2}{e_0 -1 \choose i}
$$ 
for all $i,\, 1 \leq i \leq d$.
\end{theorem}

\begin{proof}
Since $\gamma(I) \geq d-1$, we have $e_1 = \sum_{n \geq
1}\lambda(I^n/JI^{n-1})$ by Theorem \ref{d-1}. Kirby proved in \cite{k} 
that for
an $\mm$-primary ideal $I$ in a 1-dimensional Cohen-Macaulay local ring,
$e_1 \leq {e_0 \choose 2}$. Also, by \cite{sv}, if $I$ is an
$\mm$-primary ideal in a 1-dimensional Cohen-Macaulay local ring with
reduction number $r$ then $r \leq e_0(I) - 1$. By hypothesis
$\gamma(I) \geq d-1$. Therefore we can choose a reduction 
$J = (a_1, \ldots, a_d)$ of $I$
such that $a_1^*, \ldots, a_{d-1}^*$ form a regular sequence in $G(I)$
and such that $r_{\bar{a_d}}(\bar{I}) = r_J(I)$, where $``-"$ denotes 
images 
modulo $(a_1, \ldots, a_{d-1})$. Then $r = r_J(I) \leq e_0(I) - 1$. 
Also $e_1(\bar{I}) \leq {e_0(\bar{I}) \choose 2}$. 
In the first part of Theorem \ref{d} we
actually proved that if $e_1(I) = \sum_{n \geq 1}
\lambda(I^n/JI^{n-1}),$ then $\gamma(I) \geq d-1,$ $e_0(\bar{I}) = 
e_0(I)$ and $e_1(\bar{I}) = e_1(I)$. Therefore $e_1 \leq {e_0 \choose 2}$.
Now, since  $\gamma(I) \geq d-1$, Theorem \ref{d-1} gives  $ 
\sum_{n \geq 1}
\lambda(I^n/JI^{n-1}) =e_1 \leq {e_0 \choose 2}.$ Therefore 
$\lambda(I^n/JI^{n-1}) \leq {e_0 \choose 2}$ for all $n \geq 1$. By
Proposition 4.3 of \cite{hm}, $e_i = \sum_{n\geq i} {n-1 \choose i-1}
\lambda(I^n/JI^{n-1})$ for all $1 \leq i \leq d$. Therefore  for all
$i,\, 1 \leq i \leq d,$ we get
\begin{eqnarray*}  
e_i &=& \sum_{n=i}^r {n-1 \choose i-1} 
\lambda(I^n/JI^{n-1}) \\ 
    &\leq& {e_0 \choose 2} \sum_{n=i}^r
{n-1 \choose i-1}
  =  {e_0 \choose 2} {r \choose i} \leq {e_0 \choose 2} {e_0 - 1 \choose i}.
\end{eqnarray*}

\end{proof}

\vspace{3mm} 

\noindent
{\it Acknowledgment :} The authors thank R. C. Cowsik and N. V. Trung
for useful comments.

\end{document}